\documentclass[a4paper,11pt,oneside,openany]{article}

\usepackage{amssymb}
\usepackage{amsfonts}

\usepackage{graphics}
\usepackage{float}
\usepackage[dvips]{graphicx,epsfig}
\usepackage{subfigure}
\usepackage{latexsym,euscript,epic,eepic}
\newcommand{\R}{{\mathbb{R}}}
\newcommand{\E}{{\mathbb{E}}}
\newcommand{\Ind}{{\mathbb{I}}}
\newcommand{\N}{{\mathbb{N}}}

\newcommand{\WH}{W^{\mbox{\rm\tiny H}}}
\newcommand{\WZ}{W^{\mbox{\rm\tiny Z}}}
\newcommand{\WUE}{W^{\mbox{\rm\tiny HZ}}}
\newcommand{\WMVUE}{W^{\mbox{\rm\tiny opt}}}
\newcommand{\muZ}{\mu^{\mbox{\rm\tiny Z}}}
\newcommand{\muW}{\mu^{\mbox{\rm\tiny W}}}
\newcommand{\muinZ}{\mu_{i,n}^{\mbox{\rm\tiny Z}}}
\newcommand{\muinW}{\mu_{i,n}^{\mbox{\rm\tiny W}}}

\newcommand{\egloi}{\stackrel{d}{=}}
\newcommand{\egdef}{\stackrel{def}{=}}

\newcommand{\cvloi}{\stackrel{d}{\rightarrow}}

\newcommand{\CQFD}
{%
\mbox{}%
\nolinebreak%
\hfill%
\rule{2mm}{2mm}%
\medbreak%
\par%
}

\newtheorem{Theo}{Theorem}
\newtheorem{prop}{Proposition}
\newtheorem{Coro}{Corollary}
\newtheorem{Lem}{Lemma}

\begin{document}

\begin{center} {\bf{A MOVING WINDOW APPROACH FOR NONPARAMETRIC ESTIMATION OF THE CONDITIONAL TAIL INDEX }}

\vspace*{3mm}
\begin{center}
Laurent Gardes and St\'ephane Girard
\end{center}
\vspace*{5mm}
\noindent INRIA Rh\^{o}ne-Alpes, projet Mistis\\
655, avenue de l'Europe, Montbonnot\\
38334 Saint-Ismier Cedex, France. \\
\end{center}

\vspace*{6mm}
 \ \\
\noindent {\bf{Abstract}} $-$ We present a nonparametric family of estimators for the tail index of a Pareto-type distribution when covariate information is available.
Our estimators are based on a weighted sum of the log-spacings between
some selected observations. This selection is achieved through a moving
window approach on the covariate domain and a random threshold on
the variable of interest.
Asymptotic normality is proved under mild regularity conditions and illustrated for some weight functions. 
Finite sample performances are presented on a real data study. \\


\vspace*{8mm}
 \ \\
\section{Introduction}

\noindent In extreme-value statistics, one of the main problems is the estimation
of the tail index associated to a random variable $Y$.
This parameter, denoted by $\gamma$, drives the distribution tail heaviness of $Y$.
For instance, when $\gamma$ is positive, the survival function of $Y$ decreases to 
zero geometrically, and the larger $\gamma$ is, the slower is the convergence.
We refer to~\cite{EMBR} for a comprehensive treatment
of extreme-value methodology in various frameworks and to~\cite{csovih98}
for an overview of the numerous works dedicated to the estimation
of the tail index.
Here, we focus on the situation where some covariate information $x$ is
recorded simultaneously with the quantity of interest $Y$.
In the general case, the tail heaviness of $Y$ given $x$ depends on $x$,
and thus the tail index is a function $\gamma(x)$ of the covariate.  
Such situations occur for instance in climatology where one may be
interested in how climate change over years might affect extreme temperatures.
Here, the covariate is univariate (the time). Bivariate examples include
the study of extremes rainfall as a function of the geographical location.

Only a few papers address the estimation of conditional tail index.
A parametric approach is considered in~\cite{smith89} where a linear trend
is fitted to the mean of an extreme-value distribution. 
We refer to~\cite{davsmi90} for other examples of parametric models.
More recently,
Hall and Tajvidi~\cite{haltaj00} proposed to mix a non-parametric estimation
of the trend with a parametric assumption on $Y$ given $x$. 
We also refer to~\cite{beigoe03} where a kind of semi-parametric estimator is
introduced for $\gamma(\psi(\beta'x))$ where $\psi$ is a known link
function and $\beta$ is interpreted as a vector of regression coefficients.
Fully non-parametric estimators are introduced in~\cite{davram00}, 
where a local polynomial fitting of the extreme-value distribution
to the extreme observations is used. In a similar spirit, spline estimates
are fitted in~\cite{chadav05} through a penalized maximum likelihood method.
In both cases, the authors focus on univariate covariates and on the
finite sample properties of the estimators.
These results are extended in~\cite{beigoe04} where local polynomials estimates
are proposed for multivariate covariates and where their asymptotic properties
are established for very regular functions $\gamma(x)$ (at least twice
continuously differentiable). 

Similarly to these authors, we investigate how to combine nonparametric
smoothing techniques with extreme-value methods in order to obtain 
efficient estimators of $\gamma(x)$. The proposed estimator is based
on a selection, thanks to a moving window approach, of the observations
to be used in the estimator of the extreme-value
index. This estimator is a weighted sum of the rescaled log-spacings
between the selected largest observations.
This approach has several advantages.
From the theoretical point of view, very few assumptions are made
on the regularity of $\gamma(x)$ and on the nature of the covariate.
A central limit theorem is established for the proposed estimator,
without assuming that $x$ is finite dimensional.
As an example, we provide the
asymptotic rate of convergence for Lipschitzian
functions $\gamma(x)$ and multidimensional covariates $x$.
From the practical point of view, the estimator is easy to compute
since it is closed-form and thus does not require optimization procedures.

Our family of nonparametric estimators is defined in Section~\ref{defest}.
In Section~\ref{secmain}, asymptotic normality properties are established, 
and links with nonparametric regression and standard extreme-value theory (without
covariate information) are highlighted.
The choice of weights is discussed in Section~\ref{choice}. 
We first present two classical choices of weights extending Hill~\cite{hil75} and 
Zipf~\cite{krares96, schste96} estimators to the conditional case.
Next, we address the problem of obtaining minimum variance and/or unbiased 
estimators, basing on the knowledge of a second order parameter. 
The practical difficulties arising when this parameter is unknown are also
discussed. An illustration on real data is provided in Section~\ref{secreal}.
Proofs are postponed to Section~\ref{secpreuves}.

\section{Estimators of the conditional tail index}
\label{defest}

Let $E$ be a metric space associated to a metric $d$. We assume that the conditional distribution function of $Y$ given $x \in E$ is
\begin{equation}
\label{modelF}
F(y,x) = 1-y^{-1/\gamma(x)}L(y,x),
\end{equation}
where $\gamma(.)$ is an unknown positive function of the covariate $x$ and, for $x$ fixed, $L(.,x)$ is a slowly varying function, {\it{i.e.}} for $\lambda > 0$,
\[ \lim_{y \to \infty} \frac{L(\lambda y,x)}{L(y,x)} = 1. \]
Given a sample $(Y_1,x_1), \ldots, (Y_n,x_n)$ of independent observations from (\ref{modelF}), our aim is to build a point-wise estimator of the function $\gamma(.)$. More precisely, for a given $t\in E$, we want to estimate $\gamma(t)$,
focusing on the case where the design points $x_1,\ldots,x_n$ are non random.
To this end, for all $r>0$, let us denote by $B(t,r)$
the ball centered at point $t$ and with radius $r$ defined by
$$
B(t,r)=\{ x \in E, \ d(x,t) \leq r \}
$$
and let $h_{n,t}$ be a positive sequence tending to zero as $n$ goes to infinity. The proposed estimate uses a moving window approach since it is based on the response variables $Y_i's$ for which the associated covariates $x_i's$ belong to the ball $B(t,h_{n,t})$.
The proportion of such design points is thus defined by
$$
\varphi(h_{n,t})=\frac{1}{n}\sum_{i=1}^n \Ind
\{ x_i\in B(t,h_{n,t}) \}
$$
and plays an important role in this study.
It describes how the design points concentrate in the neighborhood of
$t$ when $h_{n,t}$ goes to zero, similarly to the small ball probability does,
see for instance the monograph on functional data analysis~\cite{fer06}.
Thus, the nonrandom number of observations in $[\theta,\infty) \times B(t,h_{n,t})$ is given by $m_{n,t} = n \varphi(h_{n,t})$. Let $\{ Z_i(t), \ i=1,\ldots,m_{n,t}\}$ be the response variables $Y_i's$ for which the associated covariates $x_i's$ belong to the ball $B(t,h_{n,t})$ and let $Z_{1,m_{n,t}}(t) \leq \ldots \leq Z_{m_{n,t},m_{n,t}}(t)$ be the corresponding order statistics. Our family of estimators of $\gamma(t)$ is defined by
\begin{equation}
\label{family}
{\hat{\gamma}}_n(t,W) = \sum_{i=1}^{k_{n,t}} i \log \left ( \frac{Z_{m_{n,t}-i+1,m_{n,t}}(t)}{Z_{m_{n,t}-i,m_{n,t}}(t)} \right ) W\left(i/k_{n,t},t\right)  
\left/ \sum_{i=1}^{k_{n,t}}W\left(i/k_{n,t},t\right)\right., 
\end{equation}
where $k_{n,t}$ is a sequence of integers such that $1 \leq k_{n,t} < m_{n,t}$ and $W(.,t)$ a function defined on $(0,1)$ such that $\int_0^1 W(s,t)ds \neq 0$. Thus, without loss of generality, we can assume that $\int_0^1 W(s,t)ds=1$.
Note that this family of estimators is an extension of estimators proposed
in~\cite{beidieguista02} in the situation where there is no covariate information. 
In this latter case, we also refer to~\cite{csodehmas85} for the definition
of kernel estimates based on non-increasing and non-negative functions, 
and to~\cite{soumis} for a similar work dedicated to Weibull tail-distributions.
In~\cite{vih99}, Viharos discusses the choice of the weight function to obtain universal
asymptotic normality of the corresponding weighted least-squares estimator.

We also introduce the following extended family of estimators:
\begin{equation}
\label{extendedfamily}
{\tilde{\gamma}}_n(t,\muW) = \sum_{i=1}^{k_{n,t}} i \log \left ( \frac{Z_{m_{n,t}-i+1,m_{n,t}}(t)}{Z_{m_{n,t}-i,m_{n,t}}(t)} \right ) \muinW(t)  
\left/ \sum_{i=1}^{k_{n,t}} \muinW(t) \right., 
\end{equation}
where the weights $\muinW(t)$ are 
defined by $\muinW(t) = W(i/k_{n,t},t)(1+o(1))$ uniformly in $i=1,\ldots,k_{n,t}$. 

\section{Main results}
\label{secmain}

\noindent We first give all the conditions required to obtain the asymptotic normality of our estimators. In the sequel, we fix $t \in E$ such that $\gamma(t)>0$. 

\paragraph{Assumptions on the conditional distribution.}

Let $x \in E$ be fixed. Then, model~(\ref{modelF}) is well known to be equivalent to the so-called first order condition
\begin{equation}
\label{modelU}
U(y,x)\egdef\inf\{s; F(s,x)\geq 1-1/y\} = y^{\gamma(x)} \ell(y,x), 
\end{equation}
where, for $x$ fixed, $\ell(.,x)$ is a slowly varying function. The function $U(.,x)$ is said to be regularly varying with index $\gamma(x)$. 
We refer to~\cite{BING} for a detailed account on this topic. The conditions are:

\begin{itemize}
\item[{\bf{(A.1)}}] The conditional cumulative distribution $F(.,t)$ is continuous.
\item[{\bf{(A.2)}}] There exists positive constants $c_{U}$, $z_U$ and $\alpha_{U} \leq 1$ such that for all $x \in B(t,1)$,
\[ \sup_{z \geq z_U} \left | \frac{\log U(z,x)}{\log U(z,t)}-1\right | \leq c_{U} d^{\alpha_{U}}(x,t). \]
\item[{\bf{(A.3)}}] There exists a negative function $\rho(t)$ and a rate function $b(.,t)$ satisfying $b(y,t) \to 0$ as $y \to \infty$, such that for all $\lambda \geq 1$, 
\[ \log \left ( \frac{\ell(\lambda y,t)}{\ell(y,t)} \right ) = b(y,t) \frac{1}{\rho(t)}(\lambda^{\rho(t)}-1)(1+o(1)), \]
where "o" is uniform in $\lambda\geq 1$  as $y \to \infty$.
\end{itemize}
Conditions {\bf{(A.1)}} and {\bf{(A.2)}} are regularity conditions on the conditional distribution function. 
The second-order condition {\bf{(A.3)}} on the slowly varying function is 
the cornerstone to establish the asymptotic normality of tail index estimators.
It is used in~\cite{Hausler} to prove the asymptotic normality of the
Hill estimate and in~\cite{beidiegoemat99} for one of its refinements.
The second order parameter $\rho(t)<0$ tunes the rate of convergence
of $\ell(\lambda t, x)/\ell(t,x)$ to 1. The closer $\rho(t)$ is to 0, the slower
is the convergence.  The function $b(.,t)$ is usually called the
bias function, since it drives the asymptotic behavior of most
tail index estimators.
It can be shown that necessarily, $b(.,t)$ is regularly varying with index $\rho(t)$ (see~\cite{gelhaa87}). 

\paragraph{Assumptions on the weights.}
The next assumption was first introduced in~\cite{beidieguista02} to establish 
exponential approximations for the log-spacings between extreme order
statistics. 
\begin{itemize}
\item[{\bf{(B.1)}}] The function $s\to sW(s,t)$ is absolutely continuous,
 {\it i.e.} there exists a function $u(.,t)$ defined on $(0,1)$ such that 
\begin{equation}
\label{absolutecont}
sW(s,t) = \int_0^s u(\xi,t)d\xi 
\end{equation}
with, for all $j=1,\dots,k_{n,t}$,
\begin{equation}
\label{C11}
 \left | k_{n,t} \int_{(j-1)/k_{n,t}}^{j/k_{n,t}} u(\xi,t) d\xi \right | < g \left ( \frac{j}{k_{n,t}+1},t \right ), 
\end{equation}
where $g(.,t)$ is a positive continuous function defined on $(0,1)$ and satisfying
\begin{equation}
\label{C12}
 \int_0^1 \max(1,\log(1/s))g(s,t)ds < \infty. 
\end{equation}
\item[{\bf{(B.2)}}] There exists a constant $\delta >0$ such that
$\int_0^1 |W(s,t)|^{2+\delta}ds < \infty$.
\end{itemize}

\paragraph{Assumptions on the sequences $k_{n,t}$ and $h_{n,t}$.}
We assume that $k_{n,t}$ is an intermediate
sequence, which is a classical assumption in extreme-value analysis:
\begin{itemize}
\item[{\bf{(C)}}] $n \varphi(h_{n,t})/k_{n,t} \to \infty$ and $k_{n,t} \to \infty$.
\end{itemize}
Remark that {\bf{(C)}} implies $n  \varphi(h_{n,t}) \to\infty$ {\it{i.e.}} the number of points in $[\theta,\infty) \times B(t,h_{n,t})$ goes to 
infinity as the total number of points does.

\noindent In order to simplify the notations, let
$$
b_{n,t}\egdef b \left ( \frac{n\varphi(h_{n,t})}{k_{n,t}} ,t \right )
$$
and introduce the rescaled log-spacings
\[ C_{i,n}(t)\egdef i\log \left ( \frac{Z_{m_{n,t}-i+1,m_{n,t}}(t)}{Z_{m_{n,t}-i,m_{n,t}}(t)} \right ), \ i=1,\ldots,k_{n,t}, \]
such that estimator~(\ref{family}) can be rewritten as
$$
{\hat{\gamma}}_n(t,W) = \sum_{i=1}^{k_{n,t}} C_{i,n}(t) W\left(i/k_{n,t},t\right)  
\left/ \sum_{i=1}^{k_{n,t}}W\left(i/k_{n,t},t\right)\right..
$$
Besides, in the following, each vector $\{ v_{i,n}, \ i=1,\ldots,k_{n,t} \}$ is denoted by $\{v_{i,n}\}_i$. 
Our first main result establishes the exponential regression model for $\{C_{i,n}(t)\}_i$.
 \ \\
\begin{Theo}
\label{regmodel}
Suppose {\bf{(A.1)}}, {\bf{(A.2)}}, {\bf{(A.3)}}, {\bf{(B.1)}} and {\bf{(C)}} hold. Then, the random vector $\{ C_{i,n}(t)\}_i$ has the same distribution as
\begin{eqnarray*}
\left \{ \left [ \left ( \gamma(t)+b_{n,t}  \left ( \frac{i}{k_{n,t}+1} \right )^{-\rho(t)}\right )F_i+\beta_{i,n}(t)+o_{\rm{P}} \left ( b_{n,t} \right ) \right ] (1+O_{\rm{P}}(h_{n,t}^{\alpha_U}) ) \right \}_i,
\end{eqnarray*}
uniformly in $i=1,\ldots,k_{n,t}$ with
\[ \frac{1}{k_{n,t}} \sum_{i=1}^{k_{n,t}}  W\left(i/k_{n,t},t\right)  \beta_{i,n}(t) = o_{\rm{P}} \left ( b_{n,t} \right ),  \]
and where $F_1, \ldots,F_{k_{n,t}}$ are independent standard exponential variables.
\end{Theo}
 \ \\
\noindent Similar results can be found in~\cite{DGGG} for rescaled log-spacings of Weibull-type random variables, and in~\cite{beidieguista02} in the case of Pareto-type random 
variables without covariate. We also refer to~\cite{Hillplot}
for approximations of the Hill process by sums of standard exponential random variables.
In the conditional case, {\it i.e.} when covariate information is available,
only few results exist. We refer to~\cite{FHR}, Theorem~3.5.2,
for the approximation of the
nearest neighbors distribution using the Hellinger distance and to~\cite{Gango}
for the study of their asymptotic distribution.
Our second main result establishes the asymptotic normality of our estimators. 
 \ \\
\begin{Theo}
\label{asympnor}
Suppose {\bf{(A.1)}}, {\bf{(A.2)}}, {\bf{(A.3)}}, {\bf{(B.1)}}, {\bf{(B.2)}} and {\bf{(C)}} hold. If, moreover,
\begin{equation}
\label{condsuite}
 k_{n,t}^{1/2}b_{n,t} \to \lambda(t) \in \R \mbox{ and } k_{n,t}^{1/2}h_{n,t}^{\alpha_{U}} \to 0 
\end{equation} 
then
\begin{equation}
\label{loinor}
k_{n,t}^{1/2} \left ( {\hat{\gamma}}_n(t,W) - \gamma(t) -  b_{n,t} {\cal AB}(t,W) \right ) \cvloi {\cal{N}} \left ( 0,\gamma^2(t){\cal AV}(t,W) \right ),
\end{equation}
where we have defined
$$
{\cal AB}(t,W)=\int_0^1W(s,t)s^{-\rho(t)}ds \mbox{ and }
{\cal AV}(t,W)=\int_0^1W^2(s,t) ds. 
$$
\end{Theo}
It appears that the asymptotic bias involves two parts. The first one is 
given by $b_{n,t}$ and thus depends on the original distribution itself.
The second one is given by ${\cal AB}(t,W)$. This multiplicative
factor can be made small by an appropriate choice of the weighting function $W$,
see the next section. Similarly, the variance term is inversely proportional to
$k_{n,t}$, the number of observations used to build the estimator, and the multiplicative
coefficient $\gamma^2(t) {\cal AV}(t,W)$ can also be adjusted.
When $\lambda(t)\neq 0$, the first part of condition~(\ref{condsuite}) forces
the bias to be of the same order as the standard-deviation.
The second part $k_{n,t}^{1/2}h_{n,t}^{\alpha_{U}}\to 0$ is due to
the functional nature of the tail index to estimate. It imposes to the 
fluctuations of $t\to U(.,t)$ to be negligible compared to the 
standard deviation of the estimate.

\noindent The following result establishes that the estimators of the extended family~(\ref{extendedfamily}) inherits from the asymptotic distribution of estimators in family~(\ref{family}). 
\begin{Coro}
\label{poidsequi}
Under the assumptions of Theorem~\ref{asympnor},
\begin{equation}
\label{loi2}
k_{n,t}^{1/2} \left ( {\tilde{\gamma}}_n(t,\muW) - \gamma(t) -  b_{n,t} {\cal AB}(t,W) \right ) \cvloi {\cal{N}} \left ( 0,\gamma^2(t) {\cal AV}(t,W) \right ). 
\end{equation}
\end{Coro}
We now propose a precise evaluation of the rate of
convergence obtained in Theorem~\ref{asympnor} in the particular
framework of multidimensional nonparametric regression.

\begin{Coro}
\label{corospeed}
Let $E=\R^p$ and suppose {\bf{(B.1)}}, {\bf{(B.2)}} hold. 
If, moreover, $\gamma$ is $\alpha$-Lipschitzian, 
the slowly-varying function $L$ in~(\ref{modelF}) is such that $L(y,x)=1$ for all $(y,x)\in \R_+\times\R^p$ and
\begin{equation}
\label{liminf}
\liminf_{n\to\infty} \varphi(h_{n,t})/h_{n,t}^p >0,
\end{equation}
then the convergence in distribution~(\ref{loinor}) holds 
with rate $n^{\frac{\alpha}{p+2\alpha}}\eta_n$,
where $\eta_n\to0$ arbitrarily slowly.
\end{Coro}
Condition~(\ref{liminf}) is an
assumption on the multidimensional design and on the distance $d$. Lemma~\ref{lemD}
in Section~\ref{secpreuves} provides an example of design fulfilling
this assumption.
Under the condition $L(y,x)=1$ for all $(y,x)\in \R_+\times\R^p$, 
estimating $\gamma(x)$ is a nonparametric regression problem since 
$\gamma(x)=\E(\log Y| X=x)$.
Let us highlight that the convergence rate provided by Corollary~\ref{corospeed} is,
 up to the $\eta_n$ factor, the optimal convergence
rate for estimating $\alpha$-Lipschitzian regression function
in $\R^p$, see~\cite{Stone}.

\section{Discussion on the choice of the weights}
\label{choice}

In order to illustrate the usefulness of our results, we first provide two examples
of weights extending classical extreme index estimators to the presence of
covariates.
Second, we propose some "optimal" choices of weights in the theoretical situation
where the second order parameter $\rho(t)$ is known.
Finally, we give some ideas to overcome this restrictive assumption.

\subsection{Two classical examples of weights}
\label{classical}

\noindent We first introduce an adaptation of Hill estimator to take into account
the covariate information. Considering in~(\ref{family}) the constant weight function $\WH(s,t) = 1$ for all $s \in [0,1]$ yields
\begin{equation}
\label{hill}
{\hat{\gamma}}_n(t,\WH) = \frac{1}{k_{n,t}} \sum_{i=1}^{k_{n,t}} 
i\log \left ( \frac{Z_{m_{n,t}-i+1,m_{n,t}}(t)}{Z_{m_{n,t}-i,m_{n,t}}(t)}\right) 
\end{equation}
which is formally the same expression as in~\cite{hil75}. 
Clearly, $\WH$ satisfies the assumptions {\bf{(B.1)}} and {\bf{(B.2)}} and then the asymptotic normality of ${\hat{\gamma}}_n(t,\WH)$ is a direct consequence of Theorem~\ref{asympnor}.
\begin{Coro}
\label{asymphill}
Under {\bf{(A.1)}}, {\bf{(A.2)}}, {\bf{(A.3)}}, {\bf{(C)}} and (\ref{condsuite}), the convergence in distribution~(\ref{loinor}) holds for ${\hat{\gamma}}_n(t,\WH)$ with
${\cal AB}(t,\WH)=1/(1-\rho(t))$ and  ${\cal AV}(t,\WH)=1$.
\end{Coro}

\noindent Similarly, we define a Zipf estimator (proposed simultaneously by Kratz and Resnick~\cite{krares96} and Schultze and Steinebach~\cite{schste96}) adapted to our framework. Remarking that the pairs
\[ \left ( \tau_{i,n}(t) \egdef \sum_{j=i}^{m_{n,t}} \frac{1}{j}, \log(Z_{m_{n,t}-i+1,m_{n,t}}(t))\right ), \ i=1,\ldots,m_{n,t}, \]
are approximatively distributed on a line of slope $\gamma(t)$ at least for small values of $i$ and for $h_{n,t}$ close to zero, one can propose a least-square estimator based on the $k_{n,t}$ largest observations :
\begin{equation}
\label{zipfbrut}
 {\tilde{\gamma}}_n(t,\muZ) = \sum_{i=1}^{k_{n,t}} (\tau_{i,n}(t)-{\bar{\tau}}_n(t)) \log(Z_{m_{n,t}-i+1,m_{n,t}}(t))  \Bigg / \sum_{i=1}^{k_{n,t}} (\tau_{i,n}(t)-{\bar{\tau}}_n(t)) \tau_{i,n}(t), 
\end{equation}
where $
{\bar{\tau}}_n(t) = \frac{1}{k_{n,t}} \sum_{i=1}^{k_{n,t}} \tau_{i,n}(t)$. 
Since (\ref{zipfbrut}) can be rewritten as
\[ {\tilde{\gamma}}_n(t,\muZ) = \sum_{i=1}^{k_{n,t}} i\log \left ( \frac{Z_{m_{n,t}-i+1,m_{n,t}}(t)}{Z_{m_{n,t}-i,m_{n,t}}(t)} \right ) \muinZ \left ( t \right ) \Bigg / \sum_{i=1}^{k_{n,t}} \muinZ \left ( t \right ), \]
with
\[ \muinZ(t) = \frac{1}{i} \sum_{j=1}^i (\tau_{j,n}(t)-{\bar{\tau}}_n(t)) = -\log \left ( i/k_{n,t} \right ) (1+o(1)), \]
uniformly in $i=1,\ldots,k_{n,t}$ (see Section~\ref{secpreuves} for a proof),
it appears that this estimator belongs to the extended family (\ref{extendedfamily}) associated to the weight function $\WZ(s,t)=-\log(s)$. 
Lemma~\ref{lemC} in Section~\ref{secpreuves} shows that condition {\bf{(B.1)}} is fulfilled
with $g(s,t)=1-\log(s)$ and thus Corollary~\ref{poidsequi} yields
\begin{Coro}
\label{asympzipf}
Under {\bf{(A.1)}}, {\bf{(A.2)}}, {\bf{(A.3)}}, {\bf{(C)}} and (\ref{condsuite}), the convergence in distribution~(\ref{loi2}) holds for ${\tilde{\gamma}}_n(t,\muZ)$ with
${\cal AB}(t,\WZ)=1/(1-\rho(t))^2$ and  ${\cal AV}(t,\WZ)=2$.
\end{Coro}

\subsection{Theoretical choices of weights}

In this subsection, three problems are addressed: The definition of asymptotically
unbiased estimators, of minimum variance estimators and of minimum variance 
asymptotically unbiased estimators.

\paragraph{Asymptotically unbiased estimators.} 
We propose to combine two weights functions in order to cancel the asymptotic bias. 
More precisely, we use the following result, which proof is straightforward.

\begin{prop}
\label{linearcombi}
Given two weights functions $W_1(.,t)$ and $W_2(.,t)$ satisfying {\bf{(B.1)}} and {\bf{(B.2)}} and a function $\alpha(t)$ defined on $E$, the weight function $\alpha(t)W_1(.,t)+(1-\alpha(t))W_2(.,t)$ also satisfies {\bf{(B.1)}} and {\bf{(B.2)}}.
\end{prop}

\noindent Hence, Theorem \ref{asympnor} entails that the asymptotic bias of 
the obtained estimator is given by
\[ b_{n,t} \left ( \alpha(t) {\cal AB}(t,W_1) + (1-\alpha(t))
{\cal AB}(t,W_2) \right ). \]
Clearly, if $W_1(.,t) \neq W_2(.,t)$, choosing
\begin{equation}
\label{alphaopt}
\alpha(t) = \frac{{\cal AB}(t,W_2)}{{\cal AB}(t,W_2) -{\cal AB}(t,W_1) }, 
\end{equation}
permits to cancel the asymptotic bias. As an example, one can combine the weights
of the conditional Hill and Zipf estimators defined respectively by~(\ref{hill})
and~(\ref{zipfbrut}) to obtain an asymptotically
unbiased estimator ${\hat{\gamma}}_n(t,\WUE)$ with
\[ \WUE(s,t) = \frac{1}{\rho(t)} - \left ( 1-\frac{1}{\rho(t)} \right ) \log(s). \]
The following result is a direct consequence of the above results.

\begin{Coro}
Under {\bf{(A.1)}}, {\bf{(A.2)}}, {\bf{(A.3)}}, {\bf{(C)}} and (\ref{condsuite}), the convergence in distribution~(\ref{loinor}) holds for ${\hat{\gamma}}_n(t,\WUE)$ with
${\cal AB}(t,\WUE)=0$ and  ${\cal AV}(t,\WUE)=1+(1-1/\rho(t))^2$.
\end{Coro}

\paragraph{Minimum variance estimator.}  It is also of interest to find
the weights minimizing the variance. 
The following result is the key tool to answer this question.
\begin{prop}
\label{minvar}
Let $t \in E$. The unique continuous function $W(.,t)$ such that $\int_0^1W(s,t)ds=1$ and minimizing $\int_0^1W^2(s,t)ds$ is given by $W(s,t)=1$ for all $s \in [0,1]$.
\end{prop}
It thus appears that the conditional Hill estimator~(\ref{hill}) is
the unique minimum variance estimator in~(\ref{family}).

\paragraph{Asymptotically unbiased estimator with minimum variance.} 
Finally, we provide the asymptotically unbiased estimator with minimum variance. 

\begin{prop}
\label{MVUE}
Let $t \in E$. The unique continuous function $W(.,t)$ such that
 $\int_0^1W(s,t)ds=1$, $\int_0^1W(s,t)s^{-\rho(t)}ds=0$ and minimizing $\int_0^1W^2(s,t)ds$ is given by
\[ \WMVUE(s,t) = \frac{\rho(t)-1}{\rho^2(t)} \left ( \rho(t)-1 + (1-2\rho(t)) s^{-\rho(t)}\right ). \]
\end{prop}

\noindent Remark that $\WMVUE(s,t)=\alpha(t)W_1(s,t)+(1-\alpha(t))W_2(s,t)$ with $W_1(s,t)=1$ for all $s \in (0,1)$, $W_2(s,t)=(1-\rho(t))s^{-\rho(t)}$ and $\alpha(t)=(1-\rho(t))^2/\rho^2(t)$ defined as in~(\ref{alphaopt}). 
From Lemma~\ref{lemC}, $W_1(.,t)$ and $W_2(.,t)$ both satisfy assumptions {\bf{(B.1)}} and {\bf{(B.2)}} with $g_1(s,t)=1$ and
$g_2(s,t)=(1-\rho(t))^2s^{-\rho(t)}$.
Thus, Proposition~\ref{linearcombi} and Theorem~\ref{asympnor} yield
the following corollary:
\begin{Coro}
\label{asympMVUE}
Under {\bf{(A.1)}}, {\bf{(A.2)}}, {\bf{(A.3)}}, {\bf{(C)}} and (\ref{condsuite}), the convergence in distribution~(\ref{loinor}) holds for ${\hat{\gamma}}_n(t,\WMVUE)$ with
${\cal AB}(t,\WMVUE)=0$ and  ${\cal AV}(t,\WMVUE)=(1-1/\rho(t))^2$.
\end{Coro}

\noindent Unsurprisingly, the estimators ${\hat{\gamma}}_n(t,\WUE)$ and ${\hat{\gamma}}_n(t,\WMVUE)$ requires the knowledge of the second order parameter $\rho(t)$. 
The estimation of the function $t\to\rho(t)$ is beyond the scope of this paper,
we refer to~\cite{alves1,alves2,Gomes1,Gomes2} for estimators of the second order parameter when there is no covariate information. 
The definition of estimators of the second order parameter with covariates is part of our future work as well as the study of the asymptotic properties of the $\gamma(t)$ estimator obtained by plugging the estimation of $\rho(t)$. Here, we limit ourselves to illustrating in the next subsection the effect of using a arbitrary chosen value.

\subsection{Practical choice of weights}
\label{practical}

In this subsection, we study the behavior of the estimators ${\hat{\gamma}}_{n}(t,\WUE)$ and ${\hat{\gamma}}_{n}(t,\WMVUE)$ in which we replace the second order parameter $\rho(t)$ by a arbitrary value $\rho^*<0$. We then define ${\hat{\gamma}}_n(t,\WUE_{\rho^*})$
and ${\hat{\gamma}}_n(t,\WMVUE_{\rho^*})$ with respective weights
\begin{eqnarray*}
 \WUE_{\rho^*}(s,t) &=& \frac{1}{\rho^*} - \left ( 1-\frac{1}{\rho^*} \right ) \log(s),\\
\WMVUE_{\rho^*}(s,t) &=& \frac{\rho^*-1}{(\rho^*)^2} \left ( \rho^*-1 + (1-2\rho^*) s^{-\rho^*}\right ).
\end{eqnarray*}

\noindent Their asymptotic normality is a direct consequence of Theorem~\ref{asympnor}.

\begin{figure}[H]
 \begin{center}
 {\epsfig{figure=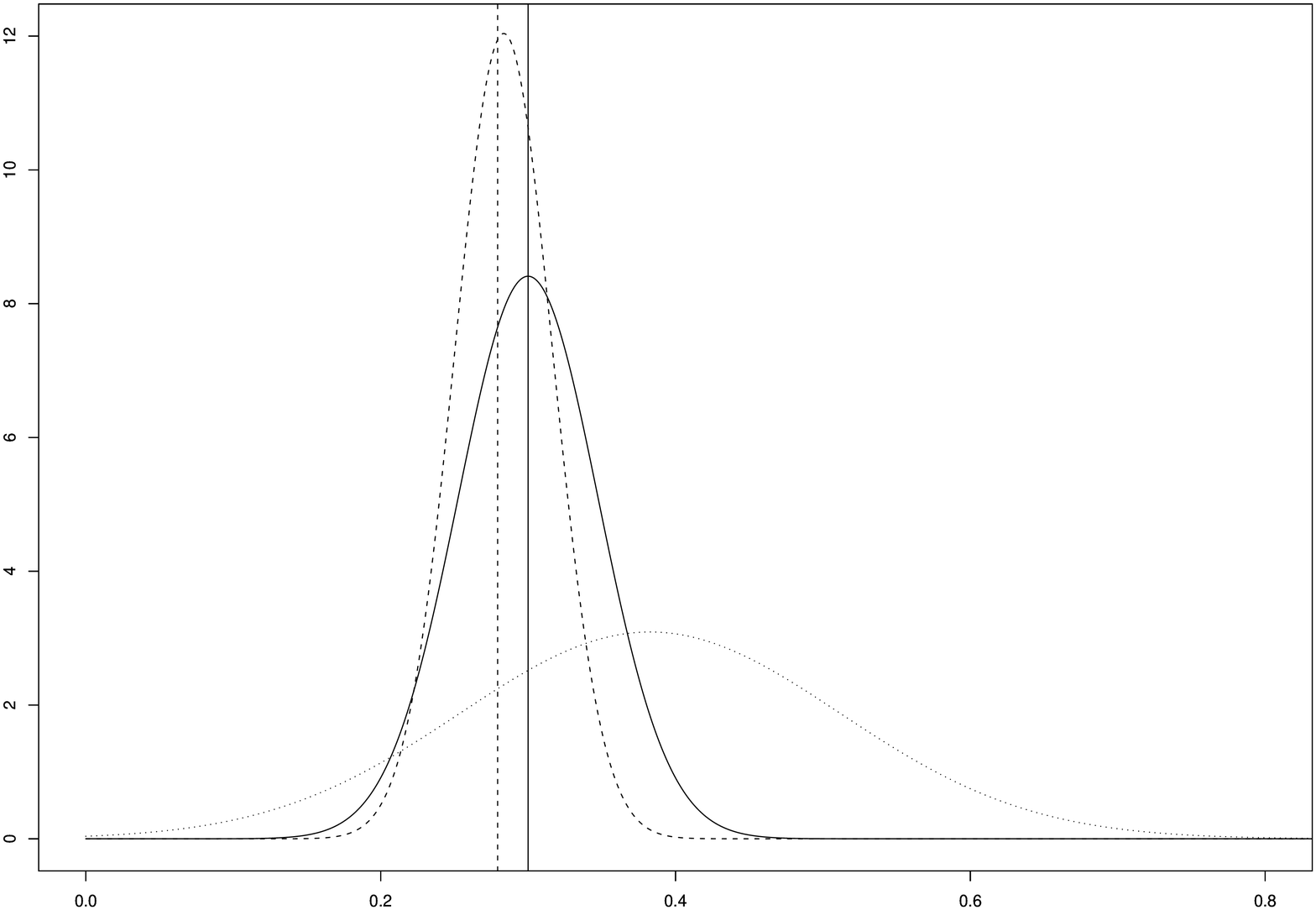,width=0.95\textwidth}}
 \caption{Densities of the asymptotic distributions of ${\hat{\gamma}}_n(t,W_{\rho^*}^{{\rm{WZ}}})$. Solid curve $\rho^*=1$, dotted curve $\rho^*=-0.2$, dashed curve $\rho^*=-5$, solid vertical line: true value $\gamma$, dotted vertical line: $\gamma + b_{n,k}{\cal{AB}}(t,W^{{\rm{Z}}})$, {\it{i.e.}}, the mean of the asymptotic distribution when $\rho^* \to -\infty$.}
 \label{figsimul}
 \end{center}
 \end{figure}

\begin{figure}[H]
 \begin{center}
 {\epsfig{figure=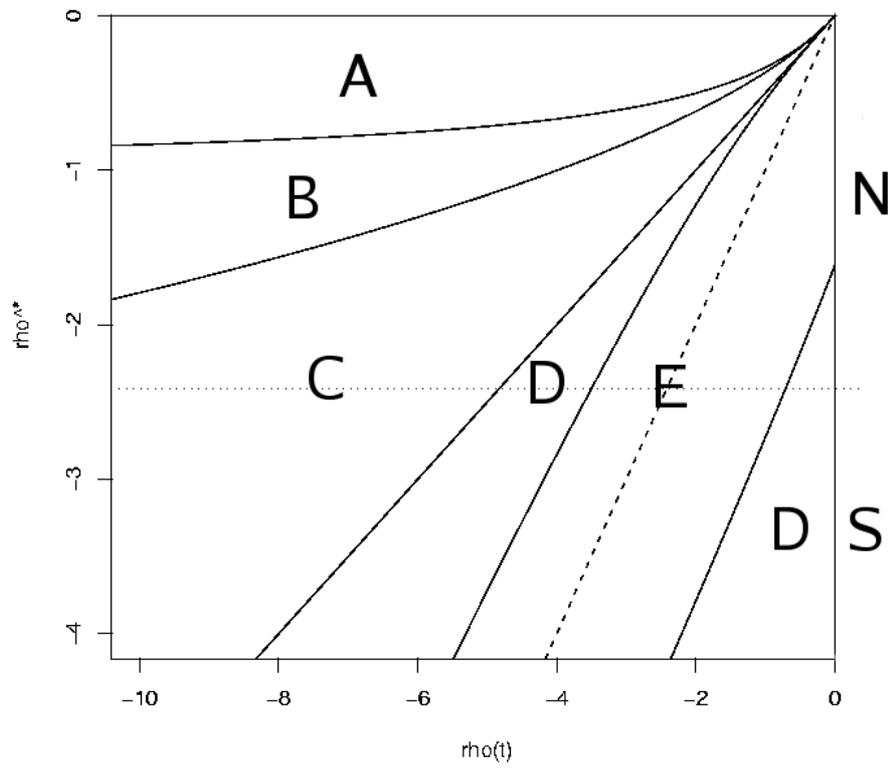,width=0.95\textwidth}}
 \caption{Comparison of the asymptotic bias and variances}
 \label{compabias}
 \end{center}
 \end{figure}

\begin{Coro}
\label{asympUE}
Under {\bf{(A.1)}}, {\bf{(A.2)}}, {\bf{(A.3)}}, {\bf{(C)}} and (\ref{condsuite}), the convergence in distribution~(\ref{loinor}) holds for ${\hat{\gamma}}_n(t,\WUE_{\rho^*})$ and  ${\hat{\gamma}}_n(t,\WMVUE_{\rho^*})$ with
$$
\begin{array}{lll}
&{\cal AB}(t,\WUE_{\rho^*})=\frac{\rho^*-\rho(t)}{\rho^*(1-\rho(t))^2}, 
&{\cal AV}(t,\WUE_{\rho^*})=1+(1-1/\rho^*(t))^2,\\
&{\cal AB}(t,\WMVUE_{\rho^*})=\frac{(1-\rho^*)(\rho^*-\rho(t))}{\rho^*(1-\rho(t))(1-\rho^*-\rho(t))}, 
&{\cal AV}(t,\WMVUE_{\rho^*})=(1-1/\rho^*(t))^2.
\end{array}
$$
\end{Coro}
The proof is a direct consequence of Theoreme~\ref{asympnor}. It appears that a bias is introduced in the asymptotic distribution. Let us also note that the asymptotic bias of the estimators ${\hat{\gamma}}_n(t,\WUE_{\rho^*})$ and ${\hat{\gamma}}_n(t,\WMVUE_{\rho^*})$ are of same sign. In term of variance, such a misspecification can allow an improvement since $\rho^* \leq \rho(t)$ yields ${\cal{AV}}(t,W^{{\rm{opt}}}_{\rho^*}) \leq {\cal{AV}}(t,W^{{\rm{opt}}})$ and ${\cal{AV}}(t,W^{{\rm{HZ}}}_{\rho^*}) \leq {\cal{AV}}(t,W^{{\rm{HZ}}})$, see Figure~\ref{figsimul}. The densities of the asymptotic distributions of ${\hat{\gamma}}_n(t,\WUE_{\rho^*})$ are represented for different choices of $\rho^*$ in case of a Burr distribution with extreme-value index $\gamma(t)=0.3$ and second order parameter $\rho(t)=-1$. Here, $m_{n,t}=5000$ and $k_{n,t}=500$ leading to $b_{n,t} \approx -0.08$. Clearly, choosing a small value of $\rho^*$ is better than choosing a large one. In fact, it is easily seen that ${\cal{AV}}(t,W^{{\rm{HZ}}}_{\rho^*}) \to {\cal{AV}}(t,W^{{\rm{Z}}})$ and ${\cal{AB}}(t,W^{{\rm{HZ}}}_{\rho^*}) \to {\cal{AB}}(t,W^{{\rm{Z}}})$ as $\rho^* \to -\infty$, whereas ${\cal{AV}}(t,W^{{\rm{HZ}}}_{\rho^*}) \to +\infty$ and ${\cal{AB}}(t,W^{{\rm{HZ}}}_{\rho^*}) \to +\infty$ as $\rho^* \to 0$. Similar conclusions hold for ${\hat{\gamma}}_n(t,W^{{\rm{opt}}})$. The consequences of the misspecification of the second order parameter on the relative efficiency are studied in~\cite{Gomes3} in the unconditional case.\\
From the practical point of view, the four estimator ${\hat{\gamma}}_n(t,W^{{\rm{H}}})$, ${\tilde{\gamma}}_n(t,\mu^{{\rm{Z}}})$, ${\hat{\gamma}}_n(t,W^{{\rm{HZ}}})$ and ${\hat{\gamma}}_n(t,W^{{\rm{opt}}})$ are easily implementable. The remainder of this paragraph is devoted to their comparison. Simple calculations lead to the following partition of the $(\rho,\rho^*)$ plane into 5 areas (see Figure~\ref{compabias}) defined as \\
\noindent A=$\{ \rho(t)<0,\rho^*<0 | \rho(t)/(2-\rho(t)) \leq \rho^*\}$, where
\[ {\cal{AB}}(t,\WZ) \leq {\cal{AB}}(t,\WH) \leq |{\cal{AB}}(t,\WUE_{\rho^*})| \leq |{\cal{AB}}(t,\WMVUE_{\rho^*})|, \]
\noindent B=$\{ \rho(t)<0,\rho^*<0 | (1-\sqrt{1-2 \rho(t)})/2 \leq \rho^* \leq \rho(t)/(2-\rho(t))\}$, where
\[ {\cal{AB}}(t,\WZ) \leq |{\cal{AB}}(t,\WUE_{\rho^*})| \leq {\cal{AB}}(t,\WH) \leq |{\cal{AB}}(t,\WMVUE_{\rho^*})|, \]
\noindent C=$\{ \rho(t)<0,\rho^*<0 | \rho(t)/2 \leq \rho^* \leq (1-\sqrt{1-2 \rho(t)})/2 \}$, where
\[ {\cal{AB}}(t,\WZ) \leq |{\cal{AB}}(t,\WUE_{\rho^*})| \leq |{\cal{AB}}(t,\WMVUE_{\rho^*})| \leq {\cal{AB}}(t,\WH), \]
\noindent D=$\{ \rho(t)<0,\rho^*<0 | \rho_1(t) \leq \rho^* \leq \rho(t)/2 \ {\rm{and}} \ \rho^* \leq \rho_2(t)\}$, where
\[ |{\cal{AB}}(t,\WUE_{\rho^*})| \leq {\cal{AB}}(t,\WZ) \leq |{\cal{AB}}(t,\WMVUE_{\rho^*})| \leq {\cal{AB}}(t,\WH), \]
\noindent E=$\{ \rho(t)<0,\rho^*<0 | \rho_2(t) \leq \rho^* \leq \rho_1(t) \}$, where
\[ |{\cal{AB}}(t,\WUE_{\rho^*})| \leq |{\cal{AB}}(t,\WMVUE_{\rho^*})| \leq {\cal{AB}}(t,\WZ) \leq {\cal{AB}}(t,\WH), \]
and with the frontier functions
\begin{eqnarray*}
 \rho_1(t)&=& \frac{\rho(t)-1-\sqrt{(1-\rho(t))^2+4(1-\rho(t))}}{2}, \\
 \rho_2(t)& =& \frac{(2+\rho(t))(\rho(t)-1)+\sqrt{(2+\rho(t))^2(1-\rho(t))^2-4\rho(t)(\rho(t)-1)(\rho(t)-2)}}{2(\rho(t)-2)}.
\end{eqnarray*}
\noindent Next, concerning the corresponding asymptotic variances, we have: \\
\noindent In the half-plane N ($\rho^* \geq -1-\sqrt{2}$),
\[ {\cal{AV}}(t,\WH) \leq {\cal{AV}}(t,\WZ) \leq {\cal{AV}}(t,\WMVUE_{\rho^*})| \leq {\cal{AV}}(t,\WUE_{\rho^*}) \]
\noindent In the half-plane S ($\rho^* \leq -1-\sqrt{2}$),
\[ {\cal{AV}}(t,\WH) \leq {\cal{AV}}(t,\WMVUE_{\rho^*}) \leq {\cal{AV}}(t,\WZ)| \leq {\cal{AV}}(t,\WUE_{\rho^*}) \]
\noindent These inequalities are summarized in Figure~\ref{compabias}.
For practical reasons, we limit $\rho(t)$ in $[-10,0]$ and $\rho^*$ in $[-4,0]$. The dashed line represents the case $\rho^*=\rho(t)$. 

\section{Illustration on real data}
\label{secreal}

In this section, we propose to illustrate our approach on 
the daily mean discharges (in cubic meters per second) of
the Chelmer river collected by the Springfield gauging station,
from 1969 to 2005. These data are provided by the
Centre for Ecology and Hydrology (United Kingdom)
and are available at {\tt http://www.ceh.ac.uk/data/nrfa}.
In this context, the variable of interest $Y$ is the daily
flow of the river and the bi-dimensional covariate $x=(x_1, x_2)$ is built
as follows: $x_1\in \{1969,1970,\dots,2005\}$ is the year
of measurement and $x_2\in\{1,2,\dots,365\}$ is the day.
The size of the dataset is $n=13,505$.

\begin{figure}[H]
 \begin{center}
 {\epsfig{figure=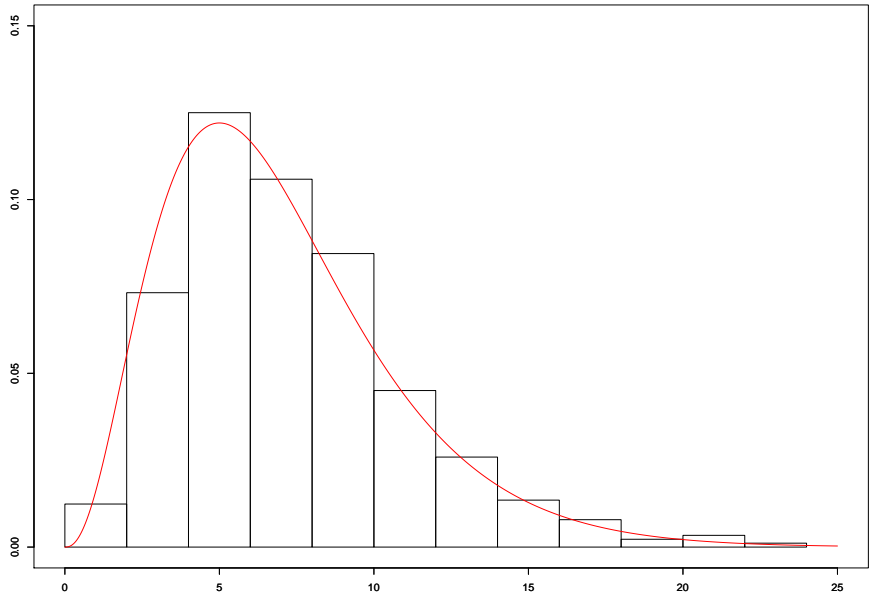,width=0.95\textwidth}}
 \caption{Histogram of the $\chi^2$ distances between the 
rescaled log-spacings and the
standard exponential distribution.
The theoretical density of the corresponding $\chi^2$ distribution
is superimposed.
}
 \label{khi2}
 \end{center}
 \end{figure}

\begin{figure}[H]
 \begin{center}
{\epsfig{figure=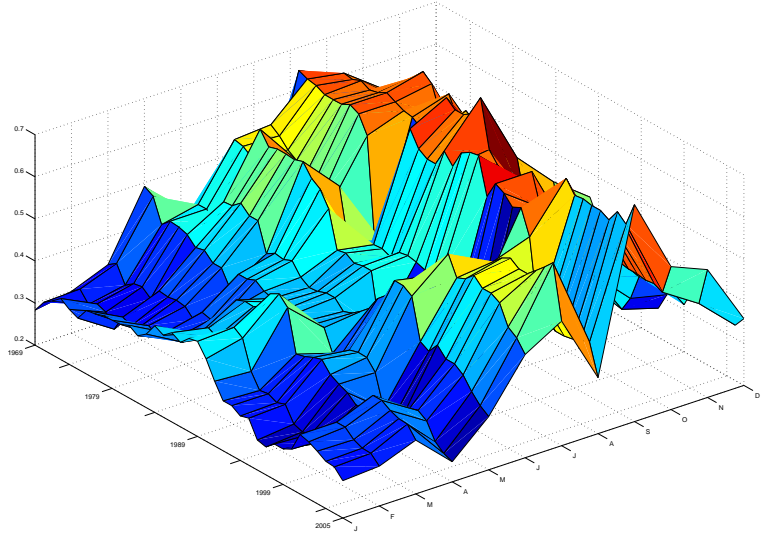,width=0.95\textwidth}}
 \caption{conditional Zipf estimator $\tilde\gamma_n(t,\muZ)$ of the tail
index computed on the
real dataset. Two covariates are available: The year ranging from 1969
to 2005 and the day ranging from 1 to 365. For the sake of readability,
only the first letter of the corresponding month is represented.
}
 \label{reszipf}
 \end{center}
 \end{figure}

The smoothing parameter $h_{n,t}$ as well as the number of upper
order statistics $k_{n,t}$ are assumed to be independent of $t$,
they are thus denoted by $h_n$ and $k_n$ respectively. They are selected 
by minimizing the following
distance between conditional Hill and Zipf estimators:
$$
\min_{h_n,k_n} \max_{t\in T} \left|\hat\gamma_n(t,\WH)-\tilde\gamma_n(t,\muZ)\right|,
$$
where $T=\{1969,1970,\dots,2005\}\times\{15,45,\dots,345\}$.
This heuristics is commonly used in functional estimation and 
relies on the idea that,
for a properly chosen pair $(h_n,k_n)$ both estimates $\hat\gamma_n(t,\WH)$
and $\tilde\gamma_n(t,\muZ)$ should yield approximately the same value.
The selected value of $h_n$ corresponds to a smoothing over 
4 years on $x_1$ and 2 months on $x_2$. Each ball $B(t,h_n)$, $t\in T$ contains
$m_n=n\varphi(h_n)=1089$ points and $k_n=54$ rescaled log-spacings are used.
This choice of $k_n$ can be validated by computing on each ball
$B(t,h_n)$, $t\in T$ the $\chi^2$ distance to the standard exponential
distribution. The histogram of these distances is superimposed
in Figure~\ref{khi2} to the theoretical density of the corresponding
$\chi^2$ distribution.
For instance, at level $5\%$, the $\chi^2$ goodness of fit test
rejects the exponential assumption in $5.7\%$ of the balls.
The resulting conditional Zipf estimator is presented on Figure~\ref{reszipf}.
The obtained values are located in the interval $[0.2, 0.7]$.
It appears that the estimated tail index is almost independent
of the year but strongly dependent of the day. The heaviest tails
are obtained in September, which means that, during this month
extreme flows are more likely than during the rest of year.

\section{Proofs}
\label{secpreuves}

\noindent For the sake of simplicity, in the sequel, we note $k_t$ for $k_{n,t}$, $b_t$ for $b_{n,t}$, $m_t$ for $m_{n,t}$ and $h_t$ for $h_{n,t}$.

\subsection{Preliminary results}

This first lemma provides sufficient conditions on $\gamma$ and $\ell$ to obtain {\bf{(A.2)}}.
\begin{Lem}
\label{lemB}
Assume that the first-order condition~(\ref{modelU}) holds. If, moreover,
there exists positive constants $z_\ell$, $c_\ell$, $c_{\gamma}$, $\alpha_{\gamma} \leq 1$ 
and $\alpha_\ell\leq 1$ such that for all $x \in B(t,1)$,
\[ |\gamma(x)-\gamma(t)| \leq c_{\gamma} d^{\alpha_{\gamma}}(x,t), \]
and
\[ \sup_{z >z_\ell} \left | \frac{\ell(z,x)}{\ell(z,t)}-1\right | \leq c_{\ell} d^{\alpha_{\ell}}(x,t), \]
then {\bf{(A.2)}} is verified with $\alpha_U=\min(\alpha_\ell,\alpha_\gamma)$.
\end{Lem}

\noindent {\bf{Proof }} $-$ Under~(\ref{modelU}), we have
$$
\frac{\log U(z,x)}{\log U(z,t)}-1 = 
\frac{(\gamma(x)-\gamma(t))\log(z) + \log\left(\frac{\ell(z,x)}{\ell(z,t)}\right)}
{\log(z)\gamma(t)\left( 1+ \frac{\log \ell(z,t)}{\gamma(t)\log(z)}\right)}.
$$
Using the well-known property of slowly varying functions 
${\log \ell(z,x) }/{\log(z)} \to 0$ as $z\to\infty$,
and taking into account that $\gamma(t)>0$, it follows that, for $z$
large enough, there exists a constant $c'_{\gamma}>0$ such that
\begin{eqnarray*}
\left| \frac{\log U(z,x)}{\log U(z,t)}-1 \right|&\leq& 
 \frac{c'_\gamma}{\gamma(t)} d^{\alpha_\gamma}(x,t) + \left|\log\left(\frac{\ell(z,x)}{\ell(z,t)}\right)\right|\\
&\leq&
 \frac{c'_\gamma}{\gamma(t)} d^{\alpha_\gamma}(x,t) + 2 \left|\frac{\ell(z,x)}{\ell(z,t)}-1\right|,
\end{eqnarray*}
since $|u| > 1/2$ entails $|\log u|\leq 2 |u-1|$. Thus,
$$
\left| \frac{\log U(z,x)}{\log U(z,t)}-1 \right|\leq
 \frac{c'_\gamma}{\gamma(t)} d^{\alpha_\gamma}(x,t) + 2 c_\ell d^{\alpha_\ell}(x,t),
$$
and the conclusion follows.\CQFD

\noindent The next lemma provides sufficient conditions on the weights to verify condition~{\bf{(B.1)}}. 
\begin{Lem}
\label{lemC}
Let $W(.,t)$ be a differentiable function on $(0,1)$. If $sW(s,t)\to0$ as $s\to 0$ then~(\ref{absolutecont}) holds with
$u(s,t)= \partial sW(s,t) / \partial s $.
Furthermore, if there exists a positive and monotone function $\phi(.,t)$ defined on $(0,1)$ such that $\max(|u(s,t)|,|W(s,t)|) \leq \phi(s,t)$, $\phi(1,t)<\infty$ and $\phi(.,t)$ is integrable at the origin then~(\ref{C11}) and~(\ref{C12}) are satisfied.
\end{Lem}

\noindent {\bf{Proof }} $-$ Clearly, since $W(.,t)$ is a differentiable function with $sW(s,t) \to 0$ as $s \to 0$, the function $sW(s,t)$ is absolutely continuous with $u(s,t)=\partial sW(s,t) / \partial s $. Furthermore, for all $j=2,\dots,k_t$, 
$$
 \left | k_t \int_{(j-1)/k_t}^{j/k_t} u(\xi,t) d\xi \right | \leq 
\sup_{s\in[(j-1)/k_t, j/k_t]} \phi(s,t).
$$
Since $\phi(.,t)$ is monotone on $(0,1)$, we have:
\[ 
 \left | k_t \int_{(j-1)/k_t}^{j/k_t} u(\xi,t) d\xi \right | \leq \left \{ \begin{array}{l l l}
\phi\left(\frac{j-1}{k_t},t\right)&\leq \phi\left(\frac{1}{2}\frac{j}{k_t+1},t\right) & {\rm{if}} \ \phi(.,t) \ {\mbox{is decreasing}}, \\
\phi\left(\frac{j}{k_t},t\right) &\leq \phi \left ( 2 \frac{j}{k_t+1},t \right )& {\rm{if}} \ \phi(.,t) \ {\mbox{is increasing}}. \end{array} \right.
\]
For $j=1$, we have
$$
 \left| k_t \int_0^{1/k_t} u(\xi,t) d\xi \right| = 
\left|W\left(\frac{1}{k_t},t\right)\right|\leq 
\phi\left(\frac{1}{k_t},t\right)\leq 
g\left(\frac{1}{k_t+1},t\right),
$$
where
\[ g(s,t) = \left \{ \begin{array}{l l}
\phi(s/2,t) & {\rm{if}} \ \phi(.,t) \ {\mbox{is decreasing}}, \\
\phi(2s,t) & {\rm{if}} \ \phi(.,t) \ {\mbox{is increasing}}. \end{array} \right.
\]
As a conclusion, condition~(\ref{C11}) is verified. From Cauchy-Schwartz inequality, to prove~(\ref{C12}), it only remains to 
verify that $\int_0^1 g(s,t)ds < +\infty$.  This is a consequence of the integrability of $\phi(.,t)$ at the origin.
\CQFD
\noindent
We now provide and example of a multidimensional design points and a distance $d$ satisfying
condition~(\ref{liminf}).  In simple words, Lemma~\ref{lemD} 
states that, if the $n$ covariates are distributed on a
 "rectangular" grid in $\R^p$,
the proportion of points in $B(t,h_{n,t})$ is asymptotically
proportional to the volume of this ball.
See~\cite{fer06}, Lemma~13.13 for a similar result in the random design setting.

\begin{Lem}
\label{lemD}
Let $E=\R^p$, $d(x,t)=\|x-t\|_\infty$ and
let $G$ be a $p$-dimensional cumulative distribution function associated to a density function $g$ such that $g(t) \neq 0$ for all $t$ in a bounded set. Assume that $G$ admits independent margins $G_1,\dots,G_p$, $n^{1/p}\in\N$, and define the lattice
${\cal L}=\{1,2,\dots,n^{1/p}\}^p\subset \N^p$. We define the multidimensional design by $\{x_\beta,\;\beta\in {\cal L}\}$ where $\beta=(\beta_1,\dots,\beta_p)\in\N^p$ is a multi-index and such that each coordinate of $x_\beta$ is given by
$$
(x_\beta)_j \egdef x_{\beta_j} \egdef G_j^{-1}\left(\frac{\beta_j-1}{n^{1/p}-1}\right), \ j=1,\ldots,p.
$$
Suppose $nh_t^p\to\infty$, then
$
\varphi(h_t)= (2h_t)^p g(t) (1+o(1)).
$
\end{Lem}
\noindent {\bf{Proof}} $-$ Using the above definitions, we have
\begin{eqnarray}
\varphi(h_t) &=& \frac{1}{n}\sum_{\beta\in {\cal L}} \Ind\{ \|x_\beta-t\|_\infty\leq h\} \nonumber\\
&=& \frac{1}{n} \sum_{\beta_1=1}^{n^{1/p}}\dots\sum_{\beta_p=1}^{n^{1/p}} \prod_{j=1}^p
\Ind\{t_j-h_t\leq x_{\beta_j}\leq t_j+h_t\} \nonumber\\
&=& \frac{1}{n} \prod_{j=1}^p \sum_{\beta_j=1}^{n^{1/p}} \Ind\{t_j-h_t\leq x_{\beta_j}\leq t_j+h_t\} \nonumber\\
&=& \frac{1}{n} \prod_{j=1}^p \sum_{\beta_j=1}^{n^{1/p}} \Ind\left\{G_j(t_j-h_t)\leq \frac{\beta_j-1}{n^{1/p}-1}\leq G_j(t_j+h_t)\right\}\nonumber\\
&=& (1-n^{-1/p})^p \prod_{j=1}^p \frac{1}{n^{1/p}-1} \sum_{\beta_j=1}^{n^{1/p}} 
Q_j\left(\frac{\beta_j-1}{n^{1/p}-1}\right),
\label{eqphi}
\end{eqnarray}
where we have introduced the indicator function
$$
Q_j(u)=\Ind\left\{G_j(t_j-h_t)\leq u \leq G_j(t_j+h_t)\right\}
$$
for $u\in[0,1]$. The above Riemann's sums can be approximated as
\begin{eqnarray*}
\frac{1}{n^{1/p}-1} \sum_{\beta_j=1}^{n^{1/p}} Q_j\left(\frac{\beta_j-1}{n^{1/p}-1}\right)
&=& \int_0^1 Q_j(u)du + O(n^{-1/p}) \\
&=& G_j(t_j+h_t)-G_j(t_j-h_t) +  O(n^{-1/p})\\
&=& 2h_t g_j(t_j) + o(h_t) +  O(n^{-1/p})\\
&=& 2h_t g_j(t_j)(1 + o(h_t)),
\end{eqnarray*}
since we assumed that $g(t)\neq 0$ and $n h_t^p\to 0$.
Replacing in~(\ref{eqphi}), the result follows.\CQFD

\subsection{Proofs of main results}

\noindent {\bf{Proof of  Theorem \ref{regmodel}}} $-$ 
Under {\bf{(A.1)}} we have 
$
\{ Z_{m_t-i+1,m_t}(t)\}_i \egloi \{ U(V_{i,m_t}^{-1},x_i)\}_i 
$
where $V_{1,m_t} \leq \ldots \leq V_{m_t,m_t}$ are the order statistics associated to the sample $V_1,\ldots,V_{m_t}$ of independent uniform variables. 
It follows that:
\begin{eqnarray*}
\{ \log(Z_{m_t-i+1,m_t}(t))\}_i
  & \egloi & \left \{ \log(U(V_{i,m_t}^{-1},t)) 
\left ( 1 + \frac{ \log(U(V_{i,m_t}^{-1},x_i))}{\log(U(V_{i,m_t}^{-1},t))}-1 \right ) \right \}_i\\
  & \egdef & \left \{ \log(U(V_{i,m_t}^{-1},t)) 
\left ( 1 + \varepsilon_{n,i} \right) \right\}_i.
\end{eqnarray*}
Now, assumption {\bf{(C)}} entails that for all $i=1,\dots,k_t$,
\[ V_{i,m_t}^{-1}\geq V_{k_t,m_t}^{-1} = (m_t/k_t)(1+o_{\rm{P}}(1)) \to \infty, \]
which implies that, for $n$ large enough, $V_{i,m_t}^{-1}\geq z_U$ for all 
$i=1,\dots,k_t$. Consequently, {\bf{(A.2)}} implies that
$$
\max_{i=1,\dots,k_t} |\varepsilon_{n,i}| \leq c_U h_t^{\alpha_U},
$$
we thus have
$
\{ \log(Z_{m_t-i+1,m_t}(t))\}_i \egloi \{ \log(U(V_{i,m_t}^{-1},t))
(1+ O_{\rm{P}} \left ( h_t^{\alpha_{U}} \right ) )\}_i.
$
The end of the proof is then a direct consequence of the following result (see~\cite{beidieguista02}, Theorem 2.1 and 2.2 for a proof):

\[ \left \{ i\log \left ( \frac{U(V_{i,m_t}^{-1},t)}{U(V_{i+1,m_t}^{-1},t)} \right ) \right \}_i = \left \{ \left ( \gamma(t) + b_t \left ( \frac{i}{k_t+1} \right )^{-\rho(t)} \right )F_i + \beta_{i,n}(t) + o_{\rm{P}}(b_t) \right \}_i, \]
where 
$\{ F_i\}_i \egdef \{i \log ({V_{i,m_t}^{-1}}/{V_{i+1,m_t}^{-1}}) \}_i$
are independent standard exponential variables and with (under {\bf{(B.1)}})
$$
\sum_{i=1}^{k_t} \left ( \frac{1}{i} \int_0^{i/k_t} u(s,t) ds \right ) \beta_{i,n}(t)  =  \frac{1}{k_t} \sum_{i=1}^{k_t} W \left ( i/k_t,t \right ) \beta_{i,n}(t) 
= o_{\rm{P}} (b_t). 
$$
\CQFD

\noindent {\bf{Proof of  Theorem \ref{asympnor}}} $-$ 
From Theorem~\ref{regmodel}, we have 
\begin{eqnarray*}
\left ( \sum_{i=1}^{k_t} W(i/k_t,t) \right ) {\hat{\gamma}}_n(t,\mu) & \egloi & \gamma(t)(1+O_{\rm{P}}(h_t^{\alpha_U})) \sum_{i=1}^{k_t}  W(i/k_t,t)F_i\\
&+&  (1+O_{\rm{P}}(h_t^{\alpha_U}))b_t \sum_{i=1}^{k_t}  W(i/k_t,t) \left ( \frac{i}{k_t+1} \right)^{-\rho(t)} F_i \\
 & + & (1+O_{\rm{P}}(h_t^{\alpha_U})) \sum_{i=1}^{k_t}  W(i/k_t,t) \beta_{i,n}(t)\\
& + & o_{\rm{P}}(b_t) \sum_{i=1}^{k_t} | W(i/k_t,t)|.
\end{eqnarray*}
Introducing
\[
T_{1,n} = \sum_{i=1}^{k_t}  W(i/k_t,t)(F_i-1), \,
T_{2,n} = \sum_{i=1}^{k_t}  W(i/k_t,t) \left ( \frac{i}{k_t+1} \right )^{-\rho(t)} (F_i-1),
\]
\[
T_{3,n} = \sum_{i=1}^{k_t}  W(i/k_t,t) \beta_{i,n}(t), \
T_{4,n}= b_t \sum_{i=1}^{k_t} W(i/k_t,t) \left ( \frac{i}{k_t+1} \right)^{-\rho(t)},\,
\]
\[ 
T_{5,n} = \sum_{i=1}^{k_t}  W(i/k_t,t),\,
T_{6,n} = \sum_{i=1}^{k_t}  |W(i/k_t,t)|, \,
T_{7,n} = \left(\sum_{i=1}^{k_t}  W^2(i/k_t,t)\right)^{1/2}, \]
we obtain the following expansion:
\begin{eqnarray}
\nonumber
\frac{T_{5,n}}{T_{7,n}}\left({\hat{\gamma}}_n(t,\mu) - \gamma(t) -
\frac{T_{4,n}}{T_{5,n}}\right)  & \egloi & 
\left(\gamma(t) \frac{T_{1,n}}{T_{7,n}} + b_t\frac{T_{2,n}}{T_{7,n}} 
+ \frac{T_{3,n}}{T_{7,n}}   \right )
(1+o_{\rm{P}}(h_t^{\alpha_U}))
\\ 
\label{expan}
&+& \left( \frac{T_{4,n}}{T_{7,n}} +\frac{T_{5,n}}{T_{7,n}}\right) o_{\rm{P}}(h_t^{\alpha_U}) + \frac{T_{6,n}}{T_{7,n}}o_{\rm{P}}(b_t).
\end{eqnarray}
Let $\delta$ be defined by {\bf (C.2}). 
From Lindeberg theorem, a sufficient condition for
$T_{1,n}/T_{7,n}\cvloi {\cal{N}}(0,1)$ is that
\begin{equation}
\label{lind}
\sum_{i=1}^{k_t}|W(i/k_t,t)|^{2+\delta} / T_{7,n}^{2+\delta} \to 0. 
\end{equation}
Since, for any integrable function $\psi$, the following convergence of Riemann sum holds,
\begin{equation}
\label{one}
 \frac{1}{k_t} \sum_{i=1}^{k_t}  \psi\left( \frac{i}{k_t}\right) \to \int_0^1 \psi(s) ds
\end{equation}
it follows that
$
T_{7,n}=k_t^{1/2} {\cal AV}(t,W)^{1/2}(1+o(1)).
$
Thus, using again~(\ref{one}),
$$
\sum_{i=1}^{k_t}|W(i/k_t,t)|^{2+\delta} / T_{7,n}^{2+\delta} 
= O(k_t^{-\delta/2}),
$$
showing that condition (\ref{lind}) is satisfied and 
\begin{equation}
\label{T1}
T_{1,n}/T_{7,n} \cvloi {\cal{N}}(0,1).
\end{equation}
Next, we focus on the term $T_{2,n}/T_{7,n}$. Remarking that this term is centered, and that its variance is finite, we can conclude that
\begin{equation}
\label{T2}
T_{2,n}/T_{7,n} = O_{\rm{P}}(1). 
\end{equation}
Theorem~\ref{regmodel} shows that
\begin{equation}
\label{T3}
T_{3,n}/T_{7,n} = o_{\rm{P}}(k_t^{1/2}b_t)=o_{\rm{P}}(1). 
\end{equation}
From repeated use of~(\ref{one}), it follows that
\begin{eqnarray}
\label{T4}
T_{4,n}/T_{5,n} & = & b_t{\cal AB}(t,W) (1+o(1))\\
T_{4,n}/T_{7,n} & = & O(k_t^{1/2}b_t)=O(1) \\
T_{5,n}/T_{7,n} & = &  k_t^{1/2} {\cal AV}(t,W)^{-1/2}(1+o(1))\\
\label{T8}
T_{6,n}/T_{7,n} & = & O(k_t^{1/2}). 
\end{eqnarray}
Replacing (\ref{T2})--(\ref{T8}) in (\ref{expan}) yields
\begin{eqnarray*}
&& k_t^{1/2} {\cal AV}(t,W)^{-1/2} 
 \left ( {\hat{\gamma}}_n(t,\mu) - \gamma(t) -  b_t {\cal AB}(t,W) \right )\\ 
&\egloi &\gamma(t) T_{1,n}/T_{7,n}
+ O( k_t^{1/2} h_t^{\alpha_U}) + o_P(1),
\end{eqnarray*}
and (\ref{T1}) gives the result.
\CQFD

\noindent {\bf {Proof of Corollary \ref{poidsequi}}} $-$
The proof consists in remarking that
\begin{eqnarray*}
 & \ & {\tilde{\gamma}}_n(t,\mu) - \gamma(t) - b_t {\cal AB}(t,W) \\
 \ & = & \sum_{i=1}^{k_t} \mu_{i,n}(t) \left ( C_{i,n}(t)-\gamma(t) - b_t{\cal AB}(t,W)  \right ) \Bigg / \sum_{i=1}^{k_t} \mu_{i,n}(t) \\
 \ & = & (1+o(1)) \sum_{i=1}^{k_t} W\left(i/k_t,t\right) \left ( C_{i,n}(t)-\gamma(t) - b_t {\cal AB}(t,W) \right) \Bigg / \sum_{i=1}^{k_t}  W\left(i/k_t,t\right)   \\
 \ & = & (1+o(1)) \left ( {\hat{\gamma}}_n(t,W) - \gamma(t) - b_t {\cal AB}(t,W)  \right ),
\end{eqnarray*}
and the conclusion follows from Theorem~\ref{asympnor}.
\CQFD

\noindent {\bf{Proof of Corollary \ref{corospeed}}} $-$ 
Assuming that $L(y,x)=1$ for all $(y,x)\in \R_+\times\R^p$ 
implies $\ell(y,x)=1$ in~(\ref{modelU})
and thus {\bf (A.3)} holds with $b(y,t)=0$. Furthermore, {\bf (A.1)}
is straightforwardly true and since $\gamma$ is $\alpha$-Lipschitzian, Lemma~\ref{lemB} entails that {\bf (A.2)} holds.
Choosing 
$h_{n,t}=n^{-\frac{1}{p+2\alpha}}$
and $k_{n,t}=n^{\frac{2\alpha}{p+2\alpha}}\eta^2_n$,
where $\eta_n\to0$ arbitrarily slowly, condition {\bf (C)}
is verified since $n h_{n,t}^p/k_{n,t}\to\infty$ and~(\ref{liminf})
imply $n\varphi(h_{n,t})/k_{n,t}\to\infty$.
As a conclusion, Theorem~\ref{asympnor}
provides the asymptotic normality of the estimator
with convergence rate $n^{\frac{\alpha}{p+2\alpha}}\eta_n$.\CQFD


\noindent {\bf {Proof of Corollary \ref{asympzipf}}} $-$
Let us first prove that~(\ref{zipfbrut}) belongs to the extended family~(\ref{extendedfamily}).  Remarking that
\[ \tau_{i,n}(t) = \sum_{j=i}^{k_t} \frac{1}{j} + \sum_{j=k_t+1}^{m_t} \frac{1}{j}, \]
estimator (\ref{zipfbrut}) can be rewritten as :
\begin{equation}
\label{zipftmp}
\sum_{i=1}^{k_t} (\tau_{i,n}(t)-{\bar{\tau}}_n(t)) \log(Z_{m_t-i+1,m_t}(t)/Z_{m_t-k_t,m_t}(t))\Bigg / \sum_{i=1}^{k_t} (\tau_{i,n}(t)-{\bar{\tau}}_n(t)) \sum_{j=i}^{k_t} \frac{1}{j}. 
\end{equation}
Next, since
\[ \log(Z_{m_t-i+1,m_t}(t)/Z_{m_t-k_t,m_t}(t)) = \sum_{j=i}^{k_t} \log(Z_{m_t-j+1,m_t}(t)/Z_{m_t-j,m_t}(t)), \]
inverting the sums in~(\ref{zipftmp}), it appears that~(\ref{zipfbrut}) belongs to family~(\ref{extendedfamily}) with
\[ \muinZ(t) = \frac{1}{i} \sum_{j=1}^{i} (\tau_{i,n}(t)-{\bar{\tau}}_n(t)). \]
Second, we prove that, uniformly in $i=1,\ldots,k_t$,
\begin{equation}
\label{W}
\muinZ(t) = -\log \left ( i/k_t \right ) (1+o(1)).
\end{equation} 
For the sake of simplicity, we introduce the following notation :
\[ S_{i,m_t} = \frac{1}{i} \sum_{j=1}^{i} \tau_{i,n}(t) = \frac{1}{i} \sum_{j=1}^{i} \sum_{l=j}^{m_t} \frac{1}{l}, \ i=1,\ldots,k_t, \]
so that $\muinZ(t) = S_{i,m_t} - S_{k_t,m_t}$. 
Furthermore, for $i=2,\ldots,k_t$,
\begin{eqnarray*}
S_{i,m_t} & = & \frac{1}{i} \sum_{j=1}^{i-1} \sum_{l=j}^{m_t} \frac{1}{l} + \frac{1}{i} \sum_{l=i}^{m_t} \frac{1}{l} \\
 \ & = & \frac{i-1}{i} S_{i-1,m_t} + \frac{1}{i} \sum_{l=i}^{m_t} \frac{1}{l} = S_{i-1,m_t} - \frac{1}{i} \left ( S_{i-1,m_t}-\sum_{l=i}^{m_t} \frac{1}{l} \right ),
\end{eqnarray*}
and remarking that
\[ S_{i-1,m_t}-\sum_{l=i}^{m_t} \frac{1}{l} = \frac{1}{i-1} \sum_{j=1}^{i-1} \sum_{l=j}^{i-1} \frac{1}{l} = 1, \]
we obtain the following recursive relation: $S_{i,m_t} = S_{i-1,m_t}-1/i$ for $i=2,\ldots,k_t$. We thus have a simplified expression of the weights:
\[ \muinZ(t) = \left \{ \begin{array}{l l} 
\sum\limits_{l=i+1}^{k_t} \frac{1}{l} & i=1,\ldots,k_t-1, \\
0 & i=k_t.
\end{array} \right . \]
We are now in position to evaluate the difference between $\muinZ(t)$ and $-\log(i/k_t)$. For $i=1,\ldots,k_t-1$,
\[ -\log \left ( i/k_t \right ) = \log \left ( \prod_{l=i+1}^{k_t} \frac{l}{l-1} \right ) = \sum_{l=i+1}^{k_t} \log \left ( 1+\frac{1}{l-1} \right ), \]
and consequently,
\begin{equation}
\label{diff}
-\log \left ( i/k_t \right ) - \muinZ(t) = \left \{ \begin{array}{l l} 
\sum\limits_{l=i+1}^{k_t} \left ( \log \left ( 1+\frac{1}{l-1} \right ) - \frac{1}{l} \right ) & i=1,\ldots,k_t-1, \\
0 & i=k_t.
\end{array} \right . 
\end{equation}
Remarking that for $l \geq 2$ the following inequality holds,
\[ 0 \leq \log \left ( 1+\frac{1}{l-1} \right ) - \frac{1}{l} \leq \frac{1}{l^2}, \]
we deduce from (\ref{diff}) that for $i=1,\ldots,k_t-1$,
\[ 0 \leq -\log \left ( i/k_t \right ) - \muinZ(t) \leq \sum_{l=i+1}^{k_t} \frac{1}{l^2}. \]
Furthermore, since
\[ \sum_{l=i+1}^{k_t} \frac{1}{l^2} \leq \int_i^{k_t} \frac{1}{x^2} dx = \frac{1}{i}-\frac{1}{k_t}, \]
we have for $i=1,\ldots,k_t-1$,
\[ 0 \leq 1-\muinZ(t)/\log(k_t/i) \leq -\frac{1}{\log(i/k_t)} \left ( \frac{1}{i}-\frac{1}{k_t} \right ). \]
Finally, since the sequence
\[ h(i)=-\frac{1}{\log(i/k_t)} \left ( \frac{1}{i}-\frac{1}{k_t} \right ), \ i \in [1,k_t[ \]
is decreasing, we have for $i=1,\ldots,k_t-1$
\[ 0 \leq 1-\muinZ(t)/\log(k_t/i) \leq \frac{1}{\log(k_t)} \left ( 1-\frac{1}{k_t} \right ), \]
proving that (\ref{W}) is true. 
The end of the proof is a consequence of Corollary~\ref{poidsequi} and 
Theorem~\ref{asympnor}. 
\CQFD

\noindent {\bf {Proof of Proposition \ref{minvar}}} $-$ For all $W$ such that $\int_0^1 W(s,t)ds=1$, we have
\[ \int_0^1 W^2(s,t)ds = 1 + \int_0^1 (W(s,t)-1)^2ds, \]
and thus minimizing $\int_0^1W^2(s,t)dt$ is equivalent to minimizing $\int_0^1(W(s,t)-1)^2ds$. Consequently, the solution of the constrained optimization problem is 
$W(.,t)=1$ almost everywhere on $[0,1]$. Since $W$ is assumed to be continuous,
the conclusion follows.
\CQFD

\noindent {\bf {Proof of Proposition \ref{MVUE}}} $-$ First, we easily check that the function $\WMVUE(.,t)$ is continuous, $\int_0^1 \WMVUE(s,t)ds=1$ and $\int_0^1 \WMVUE(s,t)s^{-\rho(t)}ds=0$. Next, remarking that for all continuous function $W(.,t)$ satisfying $\int_0^1 W(s,t)ds=1$ and $\int_0^1 W(s,t)s^{-\rho(t)}ds=0$, we have
\[ \int_0^1 W^2(s,t)ds = \left ( \frac{\rho(t)-1}{\rho(t)} \right )^2 + \int_0^1 (W(s,t)-\WMVUE(s,t))^2ds, \]
it appears that minimizing $\int_0^1W^2(s,t)ds$ is equivalent to minimizing $\int_0^1(W(s,t)-\WMVUE(s,t))^2ds$. Since $W(.,t)$ is continuous, the conclusion of the proof is straightforward.
\CQFD

\newpage

\end{document}